\documentclass[12pt,leqno]{article}
\usepackage{amsfonts,amsthm,amsmath}

\usepackage{amssymb}
\usepackage{color} 

\theoremstyle{plain}
\newtheorem{thm}{Theorem}[section]

\theoremstyle{definition}
\newtheorem{df}{Definition}[section]

\newtheorem{ex}{Example}[section]

\newcommand{\FF}{\mathbb{F}}
\newcommand{\ZZ}{\mathbb{Z}}

\newcommand{\NN}{\mathbb{N}}

\begin{document}

\title{{A generalization of the Tutte polynomials}
\footnote{This work was supported by JSPS KAKENHI(18K03217, 17K05164, 18K03388).}
}
\author{
Tsuyoshi Miezaki
\thanks{Faculty of Education, Ryukyu University, 
Okinawa, 903--0213, Japan 
},
Manabu Oura
\thanks{Graduate School of Natural Science and Technology, 
Kanazawa University,  
Ishikawa 920--1192, Japan 
},\\
Tadashi Sakuma
\thanks{Faculty of Science, 
Yamagata University,  
Yamagata 990--8560, Japan 
}, and 
Hidehiro Shinohara
\thanks{Institute for Excellence in Higher Education, 
Tohoku University,  
Miyagi 980--8576, Japan 
}
}

\date{}

\maketitle

\begin{abstract}
In this paper, we introduce the 
concept of the Tutte polynomials of genus $g$ and 
discuss some of its properties. 
We note that the Tutte polynomials of genus one 
are well-known Tutte polynomials. 
The Tutte polynomials are 
matroid invariants, and we claim that 
the Tutte polynomials of genus $g$ are also 
matroid invariants. 
The main result of this paper and the forthcoming paper declares that
the Tutte polynomials of genus $g$ are 
complete matroid invariants. 
\end{abstract}


\vspace{-0.15in}


\section{Introduction}
This is an announcement paper. 

Let $E$ be a set. A \textit{matroid $M$} on $E=E(M)$ is a pair 
$(E,\mathcal{I})$, where $\mathcal{I}$ is a non-empty 
family of subsets of $E$ with the properties 
\[
\begin{cases}
\mbox{(i)} &\mbox{ if } I\in \mathcal{I} \mbox{ and } J\subset I, \mbox{ then } 
J\in \mathcal{I}; \\
\mbox{(ii)} &\mbox{ if } I_1,I_2\in \mathcal{I} \mbox{ and }|I_1|<|I_2|, \\
&\mbox{ then there exists }
e\in I_2\setminus I_1 \\
&\mbox{ such that } I_1\cup \{e\}\in \mathcal{I}. 
\end{cases}
\]
Each element of the set $\mathcal{I}$ is called an \textit{independent set}. 
A matroid $(E, \mathcal{I})$ is \textit{isomorphic} to another matroid 
$(E', \mathcal{I}')$ 
if there is a bijection $\varphi$ from $E$ to $E'$ such that 
$\varphi(I)\in \mathcal{I}'$ holds for each member $I\in \mathcal{I}$, 
and  
$\varphi^{-1}(I')\in \mathcal{I}$ holds for each member $I'\in \mathcal{I}'$. 

It follows from the second axiom that all maximal independent sets 
in a matroid $M$ take the same cardinality, 
called the \textit{rank} of $M$. These maximal independent sets 
$\mathcal{B}(M)$ are called 
the \textit{bases} of $M$. 
The \textit{rank} $\rho(A)$ of an arbitrary subset $A$ of $E$ is the cardinality of the largest independent set contained in $A$. 


We provide examples below. 
\begin{ex}
Let $A$ be a $k\times n$ matrix over a finite field $\FF_q$. 
This offers a matroid $M$ on the set 
\[
E={\{z\in\ZZ \mid 1\leq z\leq n\}}
\] 
in which a set $I$ is independent if and only if the family of columns of $A$ whose indices belong to $I$ is linearly independent. 
Such a matroid is called a vector matroid. 
\end{ex}
\begin{ex}
Let 
\[
{E=\{z\in\ZZ \mid 1\leq z\leq n\}}
\] and $r$ a natural number with $1\leq r\leq n$. 
We define a matroid on $E$ by taking every $r$-element subset of 
$E$ to be a basis. 
This is known as the uniform matroid $U_{r,n}$. 
\end{ex}

The classification of the matroids is 
one of the most important problems in the theory of matroids. 
One tool to classify the matroids is the Tutte polynomials. 
Let $M$ be a matroid on the set $E$ having a rank function $\rho$. 
The Tutte polynomial of $M$ is defined as follows \cite{{T1},{T2},{T3}}: 
\begin{align*}
T(M)&:=T(M;x,y)\\
&:=\sum_{A\subset E}(x-1)^{\rho (E)-\rho (A)}(y-1)^{|A|-\rho (A)}. 
\end{align*}
For example, the Tutte polynomial of the uniform matroid $U_{r,n}$ is
\begin{align*}
T&(U_{r,n};x,y) \\
&=
\sum_{i=0}^r
\binom{n}{i}(x-1)^{r-i} +
\sum_{i=r+1}^n
\binom{n}{i}(y-1)^{i-r}.
\end{align*}
It is easy to demonstrate that $T(M;x,y)$ is a matroid invariant{.} 
{Two matroids are called {\it $T$-equivalent} if their Tutte polynomials
  are equivalent.}

It is well known that there exist two inequivalent matroids, 
which are $T$-equivalent \cite[p.~269]{Welsh}. 
We provide more examples below. 
Let 
\[
{E:=\{ z\in \ZZ\mid 1\leq z \leq n\}.} 
\]
Let us define the subsets $X_1$, $X_2$, $X_3$ of $E$ 
as follows: 
\[
{\begin{cases}
&X_1:=\{ z\in \ZZ\mid 1\leq z \leq r\}; \\
&X_2:=\{ z\in \ZZ\mid r+1\leq z \leq 2r\}; \\ 
&X_3:=\{ z\in \ZZ\mid {r}\leq z \leq {2r-1}\}. 
\end{cases}}
\]
Let $R_{r,n}$ denote the matroid on $E$ such that  
\[
\mathcal{B}(R_{r,n})=\mathcal{B}(U_{r, n}) 
\setminus \{ X_1, X_2\}, 
\]
$Q_{r,2n}$ denote the matroid on $E$ such that  
\[
\mathcal{B}(Q_{r,n})= 
\mathcal{B}(U_{r, n})\setminus \{ X_1, X_3 \}.
\] 
Then, for $2r \leq n$, $r\geq 3$, 
$R_{r,n}$ and $Q_{r,n}$ are matroids. 
Both matroids have {exactly} two dependent sets of size $r${{, namely}}
$\{ X_1, X_2\}$ {of $R_{r,n}$} and $\{ X_1, X_3 \}$ {of $Q_{r,n}$}. 
Therefore, if $R_{r,n}$ and $Q_{r,n}$ are isomorphic, 
there exists $\varphi$ such that 
\[
\varphi(\{ X_1, X_2\})=\{ X_1, X_3\}. 
\]
This is a contradiction since $\varphi$ is bijective, 
and we know that  
$R_{r,n}$ and $Q_{r,n}$ are non-isomorophic matroids. 

{On the other hand,}
\[
T(R_{r,n})= T(Q_{r,n}). 
\]
In fact, the difference between 
\[
T(U_{r,n})-T(R_{r,n})
\] 
and 
\[
T(U_{r,n})-T(Q_{r,n})
\] are 
{zero} since $R_{r,n}$ and $Q_{r,n}$ are obtained from 
$U_{r,n}$ after deleting the two maximal independent sets.

This gives rise to a natural question: is there
a generalization of the Tutte polynomial {which 
identifies such $T$-equivalent but inequivalent matroids?}
This paper aims to provide a candidate 
generalization that answers this. 
In Section \ref{sec:genus}, 
we provide the concept of the Tutte polynomial of genus $g$. 
In Section \ref{sec:main}, 
we provide the main results. 
The details of the proofs will be presented in our forthcoming paper \cite{MOSS}.

\section{Tutte polynomials of genus $g$}\label{sec:genus}
We now present the concept of the Tutte polynomial of genus $g$. 
\begin{df}
Let $M:=(E,\mathcal{I})$ be a matroid. 
Let 
\[
{{\Lambda_1:=\{z\in \ZZ\mid 1\leq z\leq g\}}}
\] 
and 
let 
\[
{\Lambda_2:=\displaystyle{{\Lambda_1} \choose 2}}.
\] 
For every element $\lambda \in {\Lambda_2}$, 
and $A_i\subset E$, 
let us denote 
\[
{A_{\cap(\lambda)}:=\cap_{i\in\lambda}A_i} \mbox{ and } {A_{\cup(\lambda)}:=\cup_{i\in\lambda}A_i}.
\]
Let $g$ be a natural number. 
Then, the genus $g$ of the Tutte polynomial $T^{(g)}(M)$ of the matroid $M$ will be
defined as follows:
\begin{align*}
T^{(g)}(M):=T^{(g)}(M;&x_{\lambda_1},y_{\lambda_1},x_{\cap\lambda_2},y_{\cap\lambda_2},\\
x_{\cup\lambda_2},y_{\cup\lambda_2}&
:\lambda_1\in \Lambda_1,\lambda_2\in \Lambda_2)\\
:=\sum_{A_1,\ldots,A_g \subseteq E}
\prod_{\lambda \in {\Lambda_1}}
&(x_{\lambda}-1)^{\rho(M)-\rho(A_{\lambda})}\\
&(y_{\lambda}-1)^{|A_{\lambda}|-\rho(A_{\lambda})}\\
\prod_{\lambda \in {\Lambda_2}}
&(x_{\cap(\lambda)}-1)^{\rho(M)-\rho(A_{\cap(\lambda)})}\\
&(y_{\cap(\lambda)}-1)^{|A_{\cap(\lambda)}|-\rho(A_{\cap(\lambda)})}\\
\prod_{\lambda \in {\Lambda_2}}
&(x_{\cup(\lambda)}-1)^{\rho(M)-\rho(A_{\cup(\lambda)})}\\
&(y_{\cup(\lambda)}-1)^{|A_{\cup(\lambda)}|-\rho(A_{\cup(\lambda)})}. 
\end{align*}
\end{df}
It is easy to demonstrate that $T^{(g)}(M)$ is matroid invariant and 
if two matroids have the same Tutte polynomial for genus $g$, 
we call them $T^{(g)}$-equivalent. 
For example, 
the genus $2$ for the Tutte polynomial $T^{(2)}(M)$ of the matroid $M$ is 
as follows: 
\begin{align*}
T^{(2)}(M;x_{1},&x_{2},y_{1},y_{2},x_{\cap\{12\}}
,y_{\cap\{12\}},x_{\cup\{12\}},y_{\cup\{12\}})\\
=\sum_{A_1,A_2\subset E}&
(x_1-1)^{\rho(E)-\rho(A_1)}\\
&(x_2-1)^{\rho(E)-\rho(A_2)}\\
&(y_1-1)^{|A_1|-\rho(A_1)}\\
&(y_2-1)^{|A_2|-\rho(A_2)}\\
&(x_{\cap\{1,2\}}-1)^{\rho(E)-\rho(A_1\cap A_2)}\\
&(y_{\cap\{1,2\}}-1)^{|A_1\cap A_2|-\rho(A_1\cap A_2)}\\
&(x_{\cup\{1,2\}}-1)^{\rho(E)-\rho(A_1\cup A_2)}\\
&(y_{\cup\{1,2\}}-1)^{|A_1\cup A_2|-\rho(A_1\cup A_2)}.
\end{align*}
We remark that for $g\in \NN_{\geq 2}$, 
the specialization of $T^{(g)}(M)$ is $T^{(g-1)}(M)$. 
For example, 
\[
T^{(2)}(M;x_{1},2,y_{1},2,2,2
,2)=2^{|E|}T^{(1)}(M;x_1,y_1). 
\]
Therefore, if 
\[
T^{(g)}(M)=T^{(g)}(M')
\] 
then 
\[
T^{(i)}(M)=T^{(i)}(M'), 
\]
for all $1\leq i\leq g$.


\section{Main results}\label{sec:main}
The main result of the present paper is as follows: 
\begin{thm}
\begin{enumerate}
\item 
The Tutte polynomial of genus $g$ 
{$\{T^{(g)}\}_{g=1}^{\infty}$} is a complete invariant for matroids. 
\item 
For $2r \leq n$, $r\geq 3$, 
\[
T^{(2)}(R_{r,n})\neq T^{(2)}(Q_{r,n}). 
\]

\item 

Let 
\[
E:=\{ z\in \ZZ\mid 1\leq z \leq 4n\}. 
\]
with $n\geq 3$. 
Let us define the subsets $Y_1$, $Y_2$ of $2^{E}$ 
as follows: 
\begin{align*}
Y_1&:=\{\{1,2,3\}, \{3,4,5\},\ldots, \\
&\{2n-3, 2n-2,2n-1\}, \\
&\{2n-1, 2n, 1\},\\
&\{2n+1, 2n+2, 2n+3\}, \\
&\{2n+3, 2n+4, 2n+5\},\ldots, \\
&\{4n-3, 4n-2, 4n-1\}, \\
&\{4n-1, 4n,2n+1\}\}, \\
Y_2&:=\{\{1,2,3\}, \{3,4,5\},\ldots, \\
&\{2n-3, 2n-2,2n-1\}, \\
&\{2n-1, 2n, 2n+1\},\\
&\{2n+1, 2n+2, 2n+3\}, \\
&\{2n+3, 2n+4, 2n+5\},\ldots, \\
&\{4n-3, 4n-2, 4n-1\}, \\
&\{4n-1, 4n,1\}\}.  
\end{align*}
Let $S_{4n}$ denote the independence system on $E$ such that  
\[
\mathcal{B}(S_{4n})=\mathcal{B}(U_{3, 4n}) 
\setminus  Y_1, 
\]
{$S'_{4n}$} denote the independence system on $E$ such that  
\[
\mathcal{B}({S'_{4n}})= 
\mathcal{B}(U_{3, 4n})\setminus Y_2. 
\]
Then, 
$S_{4n}$ and ${S'_{4n}}$ are matroids. 
Let 
\[
m_1=\left\lfloor\frac{-1+\sqrt{1+8n}}{2}\right\rfloor\mbox{ and } 
m_2=2\lceil\sqrt{n}\rceil. 
\]
We have 
\[
\begin{cases}
&T^{(m_1)}(S_{4n})= T^{(m_1)}({S'_{4n}});\\ 
&T^{(m_2)}(S_{4n})\neq  T^{(m_2)}({S'_{4n}}). 
\end{cases}
\]

\end{enumerate}
\end{thm}
In particular, for a matroid $M$, $T^{(|\mathcal{B}(M)|)}(M)$ 
determines $M$. For detailed explanation, see \cite{MOSS}.



\section*{Acknowledgments}
This work was supported by JSPS KAKENHI (18K03217, 17K05164, 18K03388).



\begin{thebibliography}{99}



















\bibitem{MOSS}
T. Miezaki, M.~Oura, T.~Sakuma, H.~Shinohara, 
The Tutte polynomials of genus $g$, in preparation. 







\bibitem{T1}
W.~T.~Tutte, 
A ring in graph theory, 
{\sl Proc.~Cambridge Philos.~Soc.}, 43:26--40, 1947.

\bibitem{T2}
W.~T.~Tutte, 
A contribution to the theory of chromatic polynomials, 
{\sl Canadian J.~Math.}, 6:80--91, 1954.

\bibitem{T3}
W.~T.~Tutte, 
On dichromatic polynomials,
{\sl J.~Combinatorial Theory}, 2:301--320, 1967.

\bibitem{Welsh}
D.~J.~A.~Welsh, 
Matroid Theory, 
{\sl Academic Press, London}, 1976.

\end{thebibliography}
\end{document}